# IRA assisted MMC-based topology optimization method

Kangjia Mo[a, *], Hu Wang[a, **], Zhenxing Cheng[a], Yu Li[a]

*a. State Key Laboratory of Advanced Design and Manufacturing for Vehicle Body, Hunan University, Changsha, 410082, PR China*

**Abstract** An Iterative Reanalysis Approximation (IRA) is integrated with the Moving Morphable Components (MMCs) based topology optimization (IRA-MMC) in this study. Compared with other classical topology optimization methods, the Finite Element (FE) based solver is replaced with the suggested IRA method. In this way, the expensive computational cost can be significantly saved by several nested iterations. The optimization of linearly elastic planar structures is constructed by the MMC, the specifically geometric parameters of which are taken as design variables to acquire explicitly geometric boundary. In the suggested algorithm, a hybrid optimizer based on the Method of Moving Asymptotes (MMA) approach and the Globally Convergent version of the Method of Moving Asymptotes (GCMMA) is suggested to improve convergence ratio and avoid local optimum. The proposed approach is evaluated by some classical benchmark problems in topology optimization, where the results show significant time saving without compromising accuracy.

## 1. Introduction

Topology optimization is a prominent approach to find the material distribution with the best trade-off between stiffness and volume in many fields of computational mechanics. Along with the development of topology optimization proceeding to various engineering fields, and even into commercial software, this method has faced with large-scale problems. The expensive computational cost gradually becomes a challenge for topology optimization problems. Therefore, exploring topology optimization with improved computational efficiency becomes a top priority.

[*] First author. *E-mail address*: mkj@hnu.edu.cn (K. Mo)
[**] Corresponding author. Tel.: +86 0731 88655012; fax: +86 0731 88822051.
  *E-mail address:* wanghu@hnu.edu.cn (H. Wang)

The efficiency of topology optimization depends on many factors: performance of solver, number of Degree of Freedoms (DOFs), sensitivity analysis, material modelling, etc. Up to now, extensive studies about different solvers have been applied in topology optimization to reduce the computational cost for static analysis. For example, multiple computational scales and resolutions can be employed to avoid the inherent high cost of FE analysis on a fine mesh. Nguyen et al. [1, 2] developed a Multiresolution Topology Optimization (MTOP) by employing a coarser discretization for FE and finer discretization for both density and design variables. Thomas et al. [3] proposed a new computational strategy for adaptive local mesh refinement using polygonal FE in arbitrary two-dimensional domains. Moreover, iterative solvers also play an important role in time saving. Wang et al. [4] suggested recycling parts of the search space in a Krylov subspace solver. Zakhama et al. [5] put forward a multigrid accelerated cellular automata algorithm to improve convergence of both analysis and design in topology optimization. Amir et al. [6] proposed alternative stopping criteria for a Preconditioned Conjugate Gradients (PCG) method in topology optimization so that fewer iterations were performed. Subsequently, Amir et al. [7] also exploited specific characteristics of a Multigrid Preconditioned Conjugate Gradients (MGCG) solver for reduction of cost. Besides iterative solutions, reanalysis methods are also capable of improving efficiency. Amir et al. [8-10] introduced an approximate reanalysis method, Combined Approximation (CA), which significantly reduced computational cost by at least one order of magnitude for topology optimization problems. Zheng et al. [11] integrated CA with Basic Linear Algebra Subprograms (BLAS), focusing on solutions of the eigen-value problems.

The aforementioned approaches were all implemented in the frameworks of either the Solid Isotropic Material Penalty (SIMP) technique or the Level Set Method (LSM).

Both of them have been successfully applied in widespread fields. However, there is still room for improvement due to their implicit descriptions of boundary, which makes it difficult to establish a direct link between optimization models and Computer-Aided-Design (CAD) modeling systems. The interest of CAD/CAE integration stimulated Guo et al. to propose the MMC-based topology optimization with explicit boundary descriptions in references [12-14]. Along with development, variant components were employed to extend the flexibility of geometry modeling. Different from pixel-based (e.g., SIMP) [15-18] and point-based (e.g., LSM) [19, 20] ones, the MMC-based topology optimization takes MMCs as primary building blocks. Therefore, the number of design variables, which are consist of characteristic geometry parameters of MMCs, is much smaller than conventional approaches. Despite the MMC-based topology optimization still in its infancy, it has flourished in extensive applications, such as three-dimensional problems [21, 22], geometric size control [23], additive manufacturing [24], and stress constraints [25] , etc.

In spite of the distinct superiorities of the MMC-based topology optimization, the commonplace disaster of high computational costs for static analyses cannot be circumvented either. Particularly, solving equilibrium equations will take up the dominant role of time consumption with the expansion of dimensions. It results in computational obstacles for large-scale problems. Consequently, in this study, the IRA is integrated into the MMC-based topology optimization, aiming to achieve accurate solutions efficiently. The proposed IRA consists of the exact reanalysis and the multigrid (MG) method. Reanalysis, as a kind of fast solver, is used to predict the response of modified structures efficiently without full analysis. Generally, the reanalysis methods can be classified into approximate methods and Direct Methods (DMs). CA [26-28] is one of the classical approximate reanalysis methods, which can accommodate both high-rank and local modifications, but without exact solution. DMs are usually suitable for low-rank or local modifications, and most of them can obtain exact solutions of modified structures. Theoretically, reanalysis methods only work during the relatively stable procedure of topology optimization, where the modifications are local. Thus, considering the accuracy, the proposed exact reanalysis embedded in IRA, is more appropriate for the MMC-based topology optimization. As a kind of DM, the exact reanalysis stems from the Independent Coefficients (IC) [29, 30] and Sherman-Morrison-Woodbury (SWM) formula [31, 32]. However, most of the existing reanalysis methods don't work without the direct factorization of the initial matrix. Considering large calculation and storage are required for the matrix factorization especially for large-scale problems, the MG method [33-35] is introduced. Therefore, the matrix factorization can be conducted on the coarsest grid, while relatively accurate displacement can be obtained on the finest grid. Undeniably, considerably computational saving can be promised through the IRA applied in the MMC-based topology optimization, but at the expense of oscillation near the optimal point. Sometimes the oscillation even results in the failure of convergence or converging to local minimum solutions. The diverse design variables in the MMC-based topology optimization also intensify this affect. Thereby, in order to keep robust and reliable, a hybrid optimizer combined by MMA and GCMMA replaces the single MMC or Optimality Criteria (OC) method. The series of MMA [36-43], are competent in addressing mathematical programming problems with complex objective and multiple constraints, and appropriate to solve nonlinear optimization problems, just like minimal compliance topology problems. Finally, our proposed method naturally accommodates minimum compliance problems and compliant mechanism problems, all of which can be promised with more efficient designs.

The remainder of the paper is organized as follows. In Section 2, the details of the IRA method are introduced. In Section 3, an adjoint approach for sensitivity analysis on the foundation of the MMC-based topology optimization is described and adaptive criteria about the natural implementation of the proposed algorithm framework are discussed. Sequentially, classical problems of compliance-based topology optimization are presented to demonstrate the proposed method in Section 4. Finally, conclusions are summarized in Section 5.

# 2. Iterative reanalysis approximation (IRA) method

## 2.1. Framework of the IRA method

Considering the expensive cost of analysis in topology optimization, the IRA method is suggested to improve the computational efficiency. The IRA is an integration of the multigrid (MG) and an exact reanalysis solver. The framework of the IRA method is illustrated in Fig. 1. The framework of the IRA method.

According to an iterative step of analysis in Fig. 1. The framework of the IRA method, the IRA method can be divided into two parts: the one is the two-grid level MG method in 'V-cycle'; another one is the exact reanalysis method with an effective updating strategy of Cholesky factorization for stiffness matrix in the final grid of MG. Relevant basic theories are discussed as follows.

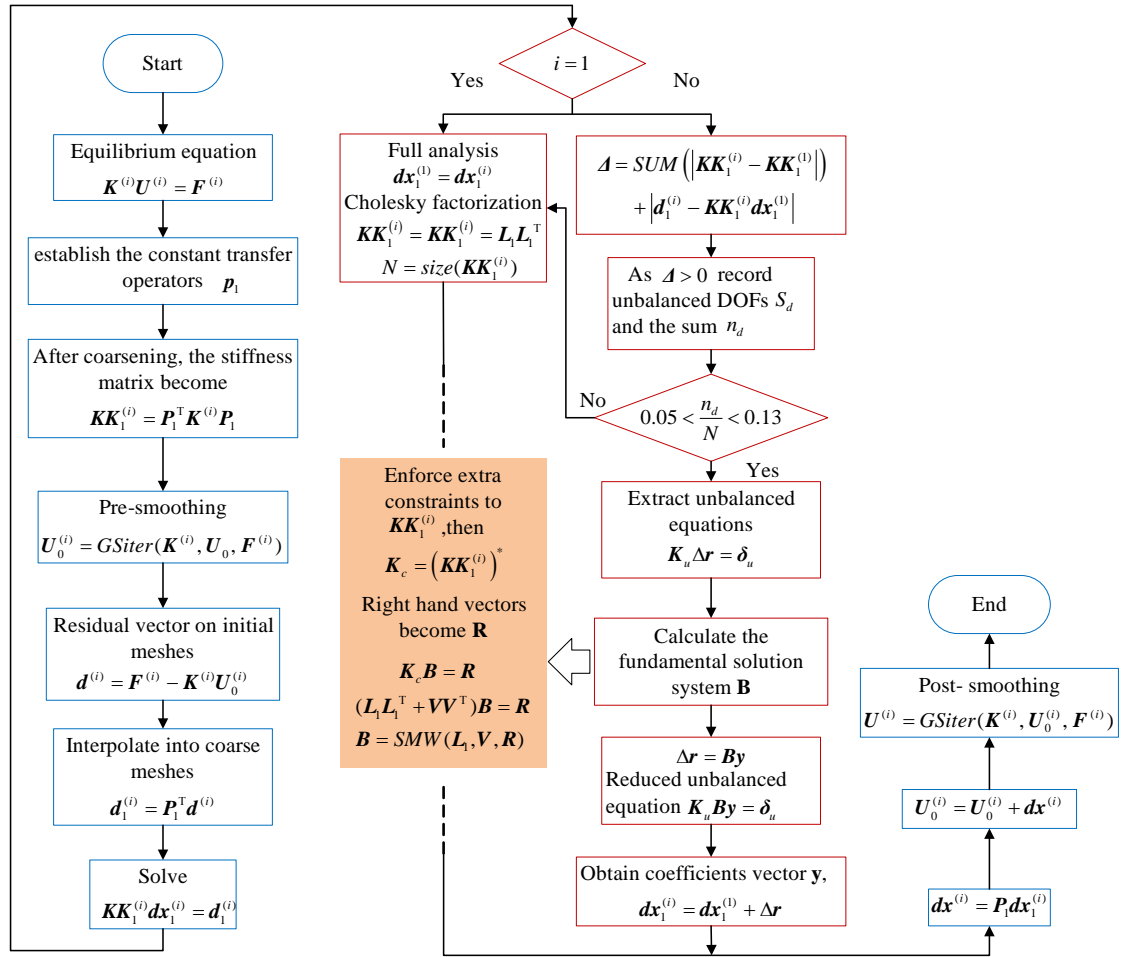

**Fig. 1.** The framework of the IRA method

## 2.2. Two-grid algorithm in 'V-cycle'

Two-grid algorithm is inherited from the MG. It can be regarded as a simplified version of the MG algorithm because only two grid levels are needed in the suggested framework. It is well known that the MG is a multilevel iterative method, which is originally proposed to solve discretized homogeneous elliptic problems. Nowadays, it has developed a family of methods that can be used to handle inhomogeneous and non-elliptic partial differential equations. Compared with ordinary iterative algorithms, the MG algorithm prevails in asymptotic convergence: the computational cost is proportional to the number of discretized variables. In contrast to other classical iterative algorithms, the MG not only has the capacity to smooth high frequency components by smoothing operator such as Gauss-Seidel or Jacobi on relatively fine grid, but also can restrict lower frequency components to a coarser grid. Ultimately, the results are obtained efficiently after only a few iterations. Meanwhile, the accuracy of solution by MG will not be deteriorated as the number of variables increases. Related numerical experiments [7] indicated that "fewer iterations are required if fewer MG levels are used, because the coarse-grid representation is more accurate". For the sake of simplicity and necessity, two-level hierarchy is suggested in the study. Owing to few levels of grid, V-cycle is more appropriate when efficiency is the dominant consideration. Generally, the MG method can be classed into startup and calculation stages. Although the stage of startup, which contains grid coarsening and establishment of transfer operators, takes relatively large computational cost. They are constructed only once and re-utilized in each cycle, since the grid is fixed in topology optimization. For the calculation stage, Gauss-Seidel smoother is embedded in the recursive iterations referred to as V-cycle, which leads to a basic grid traversal scheme.

For clarity, three primary components of the MG are stated as follows:

1. Grid coarsening

   Based on the concept of strong dependence, the partition of the coarse-grid and fine-grid variables is achieved through grid traversal scheme. Consequently, the solution domain can be divided into several grids $M_l (l=1,2,...,s)$ in the MG method, where $s$ is the number of grid level and $s=2$ in this study. For a concise recognition, the grid hierarchy with V-cycle is shown in Fig. 2.

2. Transfer operator

Transfer operator contains the prolongation $\boldsymbol{P}$ and the restriction $\boldsymbol{P}^{\mathrm{T}}$, which is the transpose of prolongation operator. The prolongation $\boldsymbol{P}$ is a $n \times n_c$ sparse matrix constructed by using standard FE interpolation on the coarse mesh for determining the weights of the fine grid nodes. Although it is time consuming to establish operators $\boldsymbol{P}$, they only need to be calculated once due to the fixed mesh in topology optimization problems. The coarse grid operator can be directly obtained by Galerkin-based coarsening according to $\boldsymbol{K}_c = \boldsymbol{P}^{\mathrm{T}} \boldsymbol{K} \boldsymbol{P}$. As a consequence, the dimensionality of equilibrium equations transfers from $n$ to $n_c$, where $n_c$ is much smaller than $n$. Herein, $n_c$ represents the number of DOFs defined on the coarse grid.

3. Smoothing operator

The process of MG starts from the classic iterative method, such as the Jacobi and the Gauss-Seidel. Nowadays, Gauss- Seidel is more commonly used on each level of MG grids except for the coarsest grid. As the smoother is being employed, it removes the high frequency components of the residual so that low order components can be effectively estimated on a coarser grid.

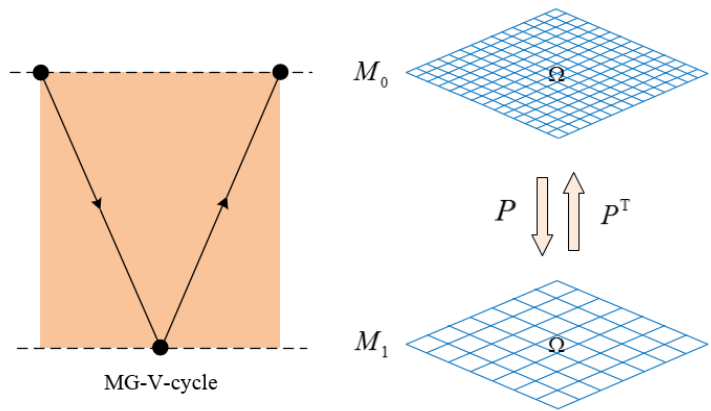

**Fig. 2.** Sketch of two grid levels of multigrid procedure

Based on the above three operators, the coarse grid correction scheme is defined as follows.

Obtain the approximate solution $U_s^{(i)}$ of $K_s^{(i)} U_s^{(i)} = F_s^{(i)}$ on the finest grid by smoother Gauss- Seidel

$$U_s^{(i)} = S^v(K_s^{(i)}, U_s^{(i)}, F_s^{(i)}), \tag{1}$$

where $S^v(K_s^{(i)}, U_s^{(i)}, F_s^{(i)})$ means returning a smoothed solution corresponding to the right hand side $F_s^{(i)}$ on the number of $s$ grid, within $v$ steps of smoothing. Superscript $i$ indicates the iterative step of the whole procedure of the IRA.

Compute the residual

$$d_s^{(i)} = F_s^{(i)} - K_s^{(i)} U_s^{(i)}. \tag{2}$$

Restrict the residual from grid $M_S$ to $M_{S+1}$ by restriction operator $P^T$

$$d_{s+1}^{(i)} = P^T d_s^{(i)}. \tag{3}$$

Compared with the standard MG, equilibrium equations on the coarsest grid are solved by the exact reanalysis method instead of the direct method, which further improves computational efficiency. The reanalysis based MG algorithm is shown in Algorithm 1.

---

**Algorithm 1**

| | |
|---|---:|
| Procedure V-cycle | |
|   For $s = 0,..., l-1$ do | Go down in the V-cycle |
|     If $s > 0$ then $x^s \leftarrow 0$ end if | |
|     Relax $K^s U^s \approx F^s$ | Apply smoother |
|     $d_{s+1} \leftarrow P_s^{s+1} d_s$ | Restrict residual |
|   End for | |
|   Solve $K^l U^l = F^l$ by the exact reanalysis method | the coarsest level |
|   For $s = l-1,..., 0$ do | Go up in the V-cycle |
|     $U^s \leftarrow U^s + P_{s+1}^s U^{s+1}$ | interpolation error & correct |
|     Relax   $U^s = U^s + S^{-1}d_s$ | Apply smoother |
|  End for | |
| End procedure | |

## 2.3. The exact reanalysis algorithm based on the coarsest level

The exact reanalysis algorithm is used to solve equilibrium equations on the coarsest level by reusing the initial solutions. Under the circumstance of local modifications in structure, the number of modified DOFs is much smaller than the unchanged DOFs. Therefore, the scale of the modified equilibrium equations solved by the reanalysis method can be remarkably reduced compared with the original ones. The details of the exact reanalysis algorithm are presented as following.

The original equilibrium equations in any iterations can be presented as

$$K^{(i)}U^{(i)} = F^{(i)}, \qquad (4)$$

where $U^{(i)}$ is the displacement in the *i-th* iterative step, $K^{(i)}$ is the stiffness matrix of system, and $F^{(i)}$ is the enforced external load. Then the equilibrium equation of the first iterative step can be represented as

$$K^{(1)}U^{(1)} = F^{(1)}. \qquad (5)$$

Although both the initial equation and the modified equation are possessed, the reanalysis algorithm cannot be directly implemented. Instead of this, two equilibrium equations should be mapped to the coarsest level firstly. After projection, Eqs.(4) and (5) can be written as

$$K_*^{(i)}dx^{(i)} = d^{(i)}, \qquad (6)$$

$$K_*^{(1)}dx^{(1)} = d^{(1)}, \qquad (7)$$

where $K_*$ is the coarsening interpolation of the MG algorithm $K_* = P^T K P$, and $d$ is obtained by Eqs.(2) and (3). Subsequently, the critical problem lies in solving the error $dx^{(i)}$ on the coarsest level.

Assume that the error $dx^{(i)}$ in the *i-th* iterative step can be decomposed into

$$dx^{(i)} = dx^{(1)} + \Delta dx \ (i > 1). \tag{8}$$

Substitute Eq.(8) into Eq.(6),

$$K_*^{(i)}(dx^{(1)} + \Delta dx) = d^{(i)} \ (i > 1), \tag{9}$$

then transform Eq.(9) into

$$K_*^{(i)}\Delta dx = d^{(i)} - K_*^{(i)}dx^{(1)} \ (i > 1). \tag{10}$$

Define the right side of Eq.(10) as

$$\delta^{(i)} = d^{(i)} - K_*^{(i)}dx^{(1)} \ (i > 1) \tag{11}$$

Thus, construct the whole solved fundament of the exact reanalysis as

$$K_*^{(i)}\Delta dx = \delta^{(i)}. \tag{12}$$

Consider only a small part of stiffness matrix will be changed when the structure is modified slightly among the contiguous iterations, thus most parts of $\delta^{(i)}$ should be zero. Thereby, Eq.(12) can be divided into balanced and unbalanced equations according to the following criterion:

$$\varDelta = sum\left(\left|K_*^{(i)} - K_*^{(1)}\right|\right) + \left|\delta^{(i)}\right|, \tag{13}$$

where $sum(\cdot)$ indicates row sum. As $\varDelta(j) > 0$, it means the $j$-th DOF is unbalanced. Therefore, the divided equations can be obtained as follows,

$$K_{*_b}^{(i)}\Delta dx = 0, \tag{14}$$

$$K_{*_u}^{(i)}\Delta dx = \delta_u^{(i)}, \tag{15}$$

where $\delta_u^{(i)}$ contains the non-zero members of $\delta^{(i)}$, and $K_{*_b}^{(i)}$ and $K_{*_u}^{(i)}$ are both some rows of $K_*^{(i)}$, but not simply upper and lower blocks.

In order to solve $\Delta dx$, the fundamental solution system of Eq. (14) can be assumed as $B$ in advance. Then the solution $\Delta dx$ can be obtained as

$$\Delta dx = By, \tag{16}$$

where $y$ is a $n_d$-dimension vector, and $n_d$ represents the number of modified DOFs.

Consequently, the critical difficulty is how to establish fundamental solution system $B$ with the help of the initial Cholesky factorization of $K_*^{(1)}$. Add extra constrains on Eq.(12) to construct the constrained linear system:

$$K_c^{(i)} B^{(i)} = R^{(i)}, \qquad (17)$$

where $K_c^{(i)}$ is converted from $K_*^{(i)}$. According to the Eq. (13), record the modified DOF $j$, then set the $j$-th row and $j$-th column of $K_*^{(i)}$ to 0, meanwhile, set the $j$-th diagonal member of $K_*^{(i)}$ to 1. The matrix $R^{(i)}$ with $n \times n_d$ is constituted by subtracting the $j$-th column of $K_*^{(i)}$, then setting the $j$-th row of it to 0, but setting the intersection (the $j$-th row corresponds to the $j$-th column of $K_*^{(i)}$ embedded in $R^{(i)}$) to 1.

According to SWM formula, Eq.(17) can be solved as:

$$B^{(i)} = \left( (L_c L_c^T)^{-1} - (L_c L_c^T)^{-1} \left( I + V(L_c L_c^T)^{-1} V^T \right)^{-1} (L_c L_c^T)^{-1} \right) R^{(i)}, \qquad (18)$$

where $L_c$ is the lower triangular matrix from restrained $L_0$, which is constrained in the same way as $K_*^{(i)}$. $L_0$ is calculated by Cholesky factorization of $K_*^{(1)}$. In this way, the factorization of the initial stiffness matrix is reused, which avoids redundant Cholesky factorization of the modified matrix iteratively. The matrix $V$ sets as following:

$$V = [I_1\ I_2\ \cdots\ I_{n_d}], \qquad (19)$$

where all the $I_1\ I_2\ \cdots\ I_{n_d}$ in Eq. (19) mean multiple rank-one modifications.

With the outcome of $B^{(i)}$, substitute Eq. (16) into Eq. (15), and obtain the reduced equations

$$K_B y = \delta_u, \qquad (20)$$

with

$$K_B = K_{*u}^{(i)} B. \qquad (21)$$

According to the reduced Eq.(20), $y$ can be directly solved. Therefore, substitute $y$

into Eq.(16), the solution of $\Delta dx$ is obtained. Owing to that the $dx^{(1)}$ is prescient, the error $dx^{(i)}$ on the coarsest level can be calculated by Eq.(8). Finally, the authentic displacement $U^{(i)}$ on the original grid can be solved by calling Algorithm 1, along the upward stage of V-cycle in the MG method.

## 3. IRA assisted MMC-based topology optimization

The integration of the IRA method with topology optimization procedure is presented in this section. Firstly, basically theoretical models associated with MMC in topology are briefly reviewed. More detailed information about these theories, can be found in [12-14, 23, 44]. Secondly, the MMC-based topology optimization with fast solver IRA is considered to formulate the minimum compliance problem and the sensitivity analysis. Finally, adaptive criteria for robust and reliable computer implementation are described.

### 3.1. MMC-based topology optimization with accurate analysis

The distinctive characteristic of the MMC-based topology optimization is that moving morphable components are taken as primary building blocks of topology optimization [12]. The structural topology descriptions about these components can be constructed in the following way:

$$\begin{cases} \phi^s(x) > 0, & \text{if } x \in \Omega^s \\ \phi^s(x) = 0, & \text{if } x \in \partial\Omega^s \\ \phi^s(x) < 0, & \text{if } x \in \mathbb{D} \setminus \Omega^s \end{cases}, \tag{22}$$

where $\mathbb{D}$ means a predefined design domain and $\Omega^s \subset \mathbb{D}$ represents a set of components $\phi^s(x) = \max(\phi_1, ..., \phi_n)$, which are made of solid material embedded in $\mathbb{D}$. The specific Topology Description Function (TDF) associated with each component can be explicitly presented with specific variables $x$ and $y$, and more details are introduced in [13]. In this study, for simplicity, the quadratic varying thickness component with straight skeleton in the form of hyperelliptic equation is adopted, which is presented in Fig. 3.

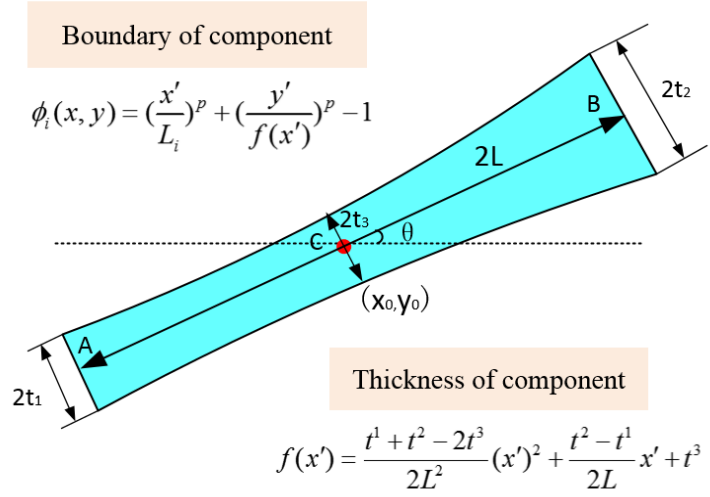

**Fig. 3**. Geometry description of the variant thickness component

The mathematic model based on the MMC under the minimum compliance objective with available volume constraint can be formulated as:

$$\left.\begin{aligned}
&Find: \boldsymbol{D} = (\boldsymbol{D}_1^T,...,\boldsymbol{D}_i^T,...,\boldsymbol{D}_{nc}^T)^T, \boldsymbol{u}(\boldsymbol{x}) \\
&Minimise: C = \int_D H(\phi^s(\boldsymbol{x};\boldsymbol{D}))\boldsymbol{f} \cdot \boldsymbol{u}dV + \int_{\Gamma_t} \boldsymbol{t} \cdot \boldsymbol{u}dS, \\
&subject\ to: \int_D (H(\phi^s(\boldsymbol{x};\boldsymbol{D})))^q \mathbb{E}:\varepsilon(\boldsymbol{u}):\varepsilon(\boldsymbol{v})dV = \int_D H(\phi^s(\boldsymbol{x};\boldsymbol{D}))\boldsymbol{f} \cdot \boldsymbol{v}dV \\
&\quad + \int_{\Gamma_t} \boldsymbol{t} \cdot \boldsymbol{v}dS, \forall \boldsymbol{v} \in u_{ad}, \\
&\quad \int_D H(\phi^s(\boldsymbol{x};\boldsymbol{D})dV \le \bar{V}, \\
&\quad \boldsymbol{D} \subset u_D, \\
&\quad \boldsymbol{u} = \bar{\boldsymbol{u}}, on\ \Gamma_u.
\end{aligned}\right\} \quad (23)$$

with

$$H = H(x) = \begin{cases} 0, & if\ x \le 0 \\ 1, & otherwise \end{cases}. \quad (24)$$

In Eq.(24), only single-phase material is considered and the Galerkin numerical method is employed. $\boldsymbol{D}_i, i=1,...,n_c$ is the vector of design variables corresponding to the $i-th$ component respectively. $H(\phi^s(\boldsymbol{x};\boldsymbol{D}))$ means utilizing Heaviside function in Eq. (24) to convert signed distance function $\phi^s(\boldsymbol{x};\boldsymbol{D})$ into solid-void (or 1-0) material distribution. Both symbols of $\boldsymbol{f}$ and $\boldsymbol{t}$ denote the body force density and the surface traction on Neumann boundary $\Gamma_t$, respectively. The parameter $q$ is an integer lager than one, and defined as $q = 2$ in this study. $\mathbb{E}$ is the constitutive matrix. The symbol $\varepsilon$ represents the second order linear strain tensor. In addition, $u$ and $v$ respectively correspond to the displacement field and the test function defined on building blocks of topology. The symbol $\bar{V}$ signifies upper bound volume of final optimized structure. $\bar{u}$ is the predefined displacement on Dirichlet boundary $\Gamma_u$ and $u_{ad} = \{v \mid v \in H^1(D), v = 0 \text{ on } \Gamma_u\}$.

**3.2. IRA-based adjoint sensitivity analysis**

For sake of simplicity, the minimum compliance topology problem with volume constraint can be concisely formulated into a general way based on work-energy principle:

$$\begin{aligned} \min \quad & C = U^T K U \\ s.t.: \quad & \int_D H(\phi^s(x; D)) dV \leq \bar{V} \quad D \subset u_D \\ \text{with}: \quad & KU = F \end{aligned} \qquad (25)$$

Taking the IRA method into consideration, the improved formulation of topology optimization problem for IRA-MMC can be rewritten as Eq. (26) at design $i-th$ iteration:

$$\begin{aligned} \min \quad & C = [U_0^{(i)} + P(dx^{(1)} + B^{(i)} y^{(i)})]^T K^{(i)} [U_0^{(i)} + P(dx^{(1)} + B^{(i)} y^{(i)})] \\ s.t.: \quad & \int_D H(\phi^s(x; D)) dV \leq \bar{V} \quad D \subset u_D \\ \text{with}: \quad & K^{(i)} [U_0^{(i)} + P(dx^{(1)} + B^{(i)} y^{(i)})] = F \\ & P^T K^{(i)} P(dx^{(1)} + B^{(i)} y^{(i)}) = P^T (F - K U_0^{(i)}) \end{aligned} \qquad (26)$$

It can be seen from Eq.(26) that the nested displacement equations are equipped with reduced system of equations by the MG and the reused initial Cholesky factorization of stiffness matrix embedded in $B^{(i)}$ simultaneously.

Within each design cycle, where the FE solver should be substituted by the IRA. Therefore, the sensitivity analysis should be calculated by the adjoint method. Introduce multipliers $\tilde{y}$ and $\tilde{\lambda}$ by Lagrange multiplier method, and then the modified objective function becomes

$$C = [U_0^{(i)} + P(dx^{(1)} + B^{(i)}y^{(i)})]^T K^{(i)}[U_0^{(i)} + P(dx^{(1)} + B^{(i)}y^{(i)})]$$
$$+ \tilde{y}\{K^{(i)}[U_0^{(i)} + P(dx^{(1)} + B^{(i)}y^{(i)})] - F\} \qquad (27)$$
$$+ \tilde{\lambda}\left[P^T K^{(i)} P(dx^{(1)} + B^{(i)}y^{(i)}) - P^T(F - KU_0^{(i)})\right]$$

In term of the minimum compliance problem, the sensitivity analysis with the nested IRA solver can be written as Eq.(28), which is differentiated with an arbitrary design variable $a$, a geometry parameter of the corresponding component.

$$\begin{aligned}
\frac{\partial C}{\partial a} &= \left[\frac{\partial U_0^{(i)}}{\partial a} + P(\frac{\partial dx^{(1)}}{\partial a} + \frac{\partial B^{(i)}}{\partial a}y^{(i)} + B^{(i)}\frac{\partial y^{(i)}}{\partial a})\right]^T K^{(i)}\left[U_0^{(i)} + P(dx^{(1)} + B^{(i)}y^{(i)})\right] \\
&+ \left[U_0^{(i)} + P(dx^{(1)} + B^{(i)}y^{(i)})\right]^T \frac{\partial K^{(i)}}{\partial a}\left[U_0^{(i)} + P(dx^{(1)} + B^{(i)}y^{(i)})\right] \\
&+ \left[U_0^{(i)} + P(dx^{(1)} + B^{(i)}y^{(i)})\right]^T K^{(i)}\left[\frac{\partial U_0^{(i)}}{\partial a} + P(\frac{\partial dx^{(1)}}{\partial a} + \frac{\partial B^{(i)}}{\partial a}y^{(i)} + B^{(i)}\frac{\partial y^{(i)}}{\partial a})\right] \\
&+ \tilde{y}K^{(i)}\left[\frac{\partial U_0^{(i)}}{\partial a} + P(\frac{\partial dx^{(1)}}{\partial a} + \frac{\partial B^{(i)}}{\partial a}y^{(i)} + B^{(i)}\frac{\partial y^{(i)}}{\partial a})\right] \\
&+ \tilde{y}\frac{\partial K^{(i)}}{\partial a}\left[U_0^{(i)} + P(dx^{(1)} + B^{(i)}y^{(i)})\right] + \tilde{\lambda}P^T \frac{\partial K^{(i)}}{\partial a}P(dx^{(1)} + B^{(i)}y^{(i)}) \\
&+ \tilde{\lambda}P^T K^{(i)} P(\frac{\partial dx^{(1)}}{\partial a} + \frac{\partial B^{(i)}}{\partial a}y^{(i)} + B^{(i)}\frac{\partial y^{(i)}}{\partial a}) + \tilde{\lambda}P^T(\frac{\partial K^{(i)}}{\partial a}U_0^{(i)} + K^{(i)}\frac{\partial U_0^{(i)}}{\partial a})
\end{aligned}, \quad (28)$$

where the transfer operator $P$ is a constant matrix, which is unnecessary to be differentiated, and $\frac{\partial dx^{(1)}}{\partial a}$ is equivalent to zero because $dx^{(1)}$ indicates the error on the coarsest grid corresponding to the first iteration.

Consider the conciseness for the explicit derivatives, only the derivative of the stiffness matrix is reserved. Consequently, the multipliers $\tilde{y}$ and $\tilde{\lambda}$ should be chosen carefully to eliminate the unnecessary items. Then the objective sensitivity analysis can be resettled as:

$$\frac{\partial C}{\partial a} = \left[ U_0^{(i)} + P(dx^{(1)} + B^{(i)} y^{(i)}) \right]^T \frac{\partial K^{(i)}}{\partial a} \left[ U_0^{(i)} + P(dx^{(1)} + B^{(i)} y^{(i)}) \right]$$

$$+ (\tilde{y} + \tilde{\lambda} P^T) \frac{\partial K^{(i)}}{\partial a} \left[ U_0^{(i)} + P(dx^{(1)} + B^{(i)} y^{(i)}) \right]$$

$$+ \left( \frac{\partial U_0^{(i)}}{\partial a} \right)^T K^{(i)} \left[ U_0^{(i)} + P(dx^{(1)} + B^{(i)} y^{(i)}) \right]$$

$$+ \left\{ \left[ U_0^{(i)} + P(dx^{(1)} + B^{(i)} y^{(i)}) \right]^T + \tilde{y} + \tilde{\lambda} P^T \right\} K^{(i)} \frac{\partial U_0^{(i)}}{\partial a} \qquad (29)$$

$$+ \left( P \frac{\partial B^{(i)}}{\partial a} y^{(i)} \right)^T K^{(i)} \left[ U_0^{(i)} + P(dx^{(1)} + B^{(i)} y^{(i)}) \right]$$

$$+ \left\{ \left[ U_0^{(i)} + P(dx^{(1)} + B^{(i)} y^{(i)}) \right]^T + \tilde{y} + \tilde{\lambda} P^T \right\} K^{(i)} P \frac{\partial B^{(i)}}{\partial a} y^{(i)}$$

$$+ \left( P B^{(i)} \frac{\partial y^{(i)}}{\partial a} \right)^T K^{(i)} \left[ U_0^{(i)} + P(dx^{(1)} + B^{(i)} y^{(i)}) \right]$$

$$+ \left\{ \left[ U_0^{(i)} + P(dx^{(1)} + B^{(i)} y^{(i)}) \right]^T + \tilde{y} + \tilde{\lambda} P^T \right\} K^{(i)} P B^{(i)} \frac{\partial y^{(i)}}{\partial a}$$

Except for the items related to $\frac{\partial K}{\partial a}$, the others in Eq. (29) need to be zero. Hence, the rest can be stacked up as:

$$\left\{ 2 \left[ U_0^{(i)} + P(dx^{(1)} + B^{(i)} y^{(i)}) \right]^T + (\tilde{y} + \tilde{\lambda} P^T) \right\} K^{(i)} \frac{\partial \left[ U_0^{(i)} + P(dx^{(1)} + B^{(i)} y^{(i)}) \right]}{\partial a} = 0. (30)$$

It is prone to be deduced from Eq. (30) that

$$\tilde{y} + \tilde{\lambda} P^T = -2 \left[ U_0^{(i)} + P(dx^{(1)} + B^{(i)} y^{(i)}) \right]^T. \qquad (31)$$

Therefore, the final expression for the objective sensitivity is

$$\frac{\partial C}{\partial a} = -\left[ U_0^{(i)} + P(dx^{(1)} + B^{(i)} y^{(i)}) \right]^T \frac{\partial K^{(i)}}{\partial a} \left[ U_0^{(i)} + P(dx^{(1)} + B^{(i)} y^{(i)}) \right]. \qquad (32)$$

This expressive form is almost identical to the one obtained by full analysis as:

$$\frac{\partial C}{\partial a} = -U^T \frac{\partial K}{\partial a} U. \qquad (33)$$

Due to the ersatz material model adopted in the MMC-based topology optimization method, the Young's modulus of this element can be interpolated as:

$$E^e = \frac{E(\sum_{i=1}^{4} (H(\phi_i^e))^q)}{4}. \qquad (34)$$

As a consequence, the complete sensitivity analysis of objective function is expressed

as:

$$\frac{\partial C}{\partial a} = -\left[\boldsymbol{U}_0^{(i)} + \boldsymbol{P}(d\boldsymbol{x}^{(1)} + \boldsymbol{B}^{(i)}\boldsymbol{y}^{(i)})\right]^T \left(\frac{E}{4}(\sum_{e=1}^{NE}\sum_{i=1}^{4} q(H(\phi_i^e))^{q-1} \frac{\partial H(\phi_i^e)}{\partial a})\boldsymbol{k}^S\right) \\ \times \left[\boldsymbol{U}_0^{(i)} + \boldsymbol{P}(d\boldsymbol{x}^{(1)} + \boldsymbol{B}^{(i)}\boldsymbol{y}^{(i)})\right] \quad , \quad (35)$$

where the $\boldsymbol{k}^S$ is the element stiffness matrix corresponding to $\phi_i^e = 1, i = 1,...,4$ and $E = 1$. The superscript $NE$ denotes the total number of elements in the ground structure. The exponent $q > 1$ is an integer, and in this study it is defined as 2.

In addition, the derivative of constrained function is shown as

$$\frac{\partial V}{\partial a} = \frac{1}{4}\sum_{e=1}^{NE}\sum_{i=1}^{4} \frac{\partial H(\phi_i^e)}{\partial a} . \quad (36)$$

### 3.3. Adaptive criteria for numerical implementation

In order to integrate the MMC-based topology optimization with the IRA seamlessly, a matrix updating strategy is suggested to activate the reanalysis procedure to keep enough accuracy in optimization. Except for the intrinsic character of IRA in improving efficiency, an additional slack convergence criterion for two-grid algorithm is adopted to strengthen this advantage. Furthermore, the optimizer is no longer the single MMA, but a hybrid optimizer constructed by MMA and GCMMA. In this way, the burdensome oscillation and instability during the stage of optimization can be alleviated. The flow chart of the IRA-MMC in the realistic implementation is demonstrated in Fig. 4.

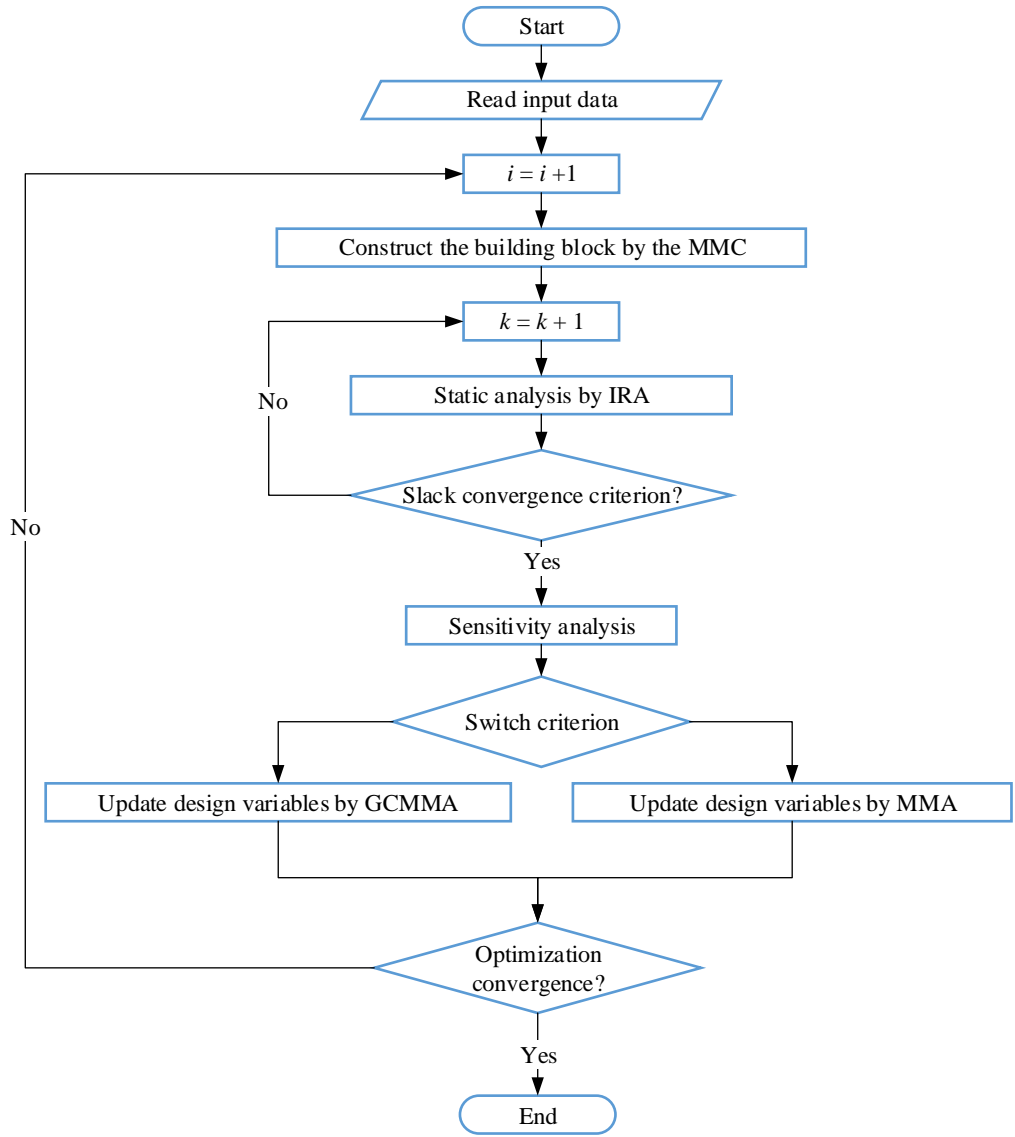

**Fig. 4.** The flowchart of the IRA-MMC based topology optimization

### 3.3.1. Strategies of matrix factorization

The primary purpose of the IRA-MMC based topology optimization is to improve efficiency with accurate results. The enough accuracy can be promised by the exact reanalysis method and MG. Although the redundant operators have been solved by MG, the implementation of the IRA method also impacts efficiency profoundly. Therefore, it is worth considering when to activate the reanalysis part of the IRA and how often to conduct a new matrix factorization. On account of dramatic structural modifications during the initial iterations in the MMC-based topology optimization,

which can be referred to

Table **1**, reanalysis can't always be utilized in all the iterations. Herein, the threshold $\eta$ for controlling matrix factorization of the IRA is described as:

$$\eta = \frac{n_d}{N} \times 100\% \tag{37}$$

The $n_d$ is the number of modified DOFs in stiffness matrix and the $N$ means the number of whole DOFs in stiffness matrix, both of which are implemented on the coarsest grid level.

Table 1 Optimization layouts of the first five iterations with respect to three examples: cantilever beam, L-shape beam and compliant mechanism

| Iterations | Cantilever beam | L-shape beam | Compliant |
|---|---|---|---|
| 1 | 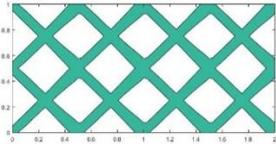 | 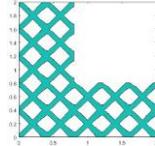 | 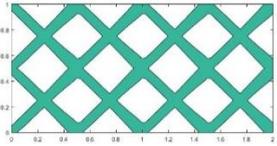 |
| 2 | 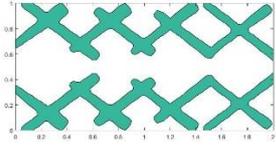 | 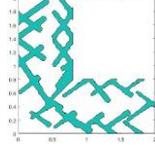 | 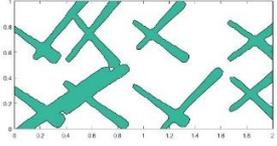 |
| 3 | 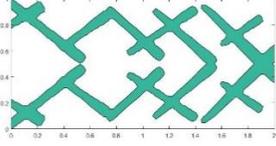 | 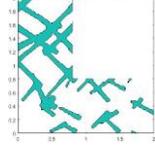 | 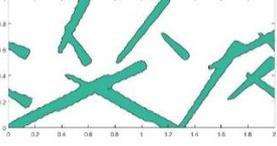 |
| 4 | 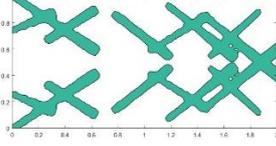 | 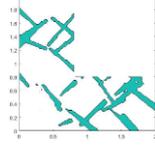 | 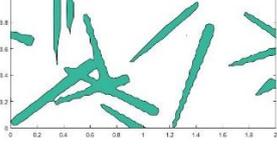 |
| 5 | 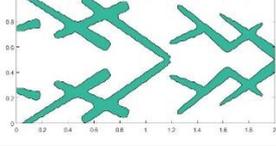 | 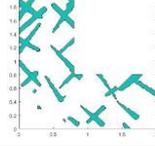 | 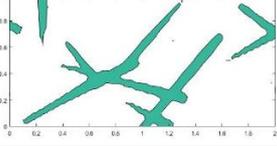 |

According to our investigation based on numerical experiments, $\eta$ is suggested to be 5%~13%. It needs to be changed carefully near the suggested values according to different problems.

### 3.3.2. Slack convergence criterion in MG

During a MG procedure, the convergence criterion generally depends on the norm of the residual $d_s^{(i)}$ in Eq.(2), which is calculated relatively to the norm of external force $F_s^{(i)}$. Once the standard criterion in Eq. (38) is satisfied, the displacement $U^{(i)}$ will converge.

$$\frac{\left\| F_s^{(i)} - K_s^{(i)} U_s^{(i)} \right\|_2}{\left\| F_s^{(i)} \right\|_2} < \varepsilon, \tag{38}$$

where the default value for the $\varepsilon$ is $10^{-6}$. However, such criterion is actually unduly conservative. Despite the relative norm of the residuals is several order of magnitude larger than $\varepsilon$, the approximate $U^{(i)}$ can be accurate enough for optimization, as shown in [6]. Consequently, it is suggested to suggest a convergence criterion based on design sensitivities as Eq.(39). This criterion not only reduces iterations in the nested IRA but also maintains sufficient accuracy by considering sensitivity analysis.

$$\frac{\left| \tilde{U}_k^{(i)T} \frac{\partial K^{(i)}}{\partial a} \tilde{U}_k^{(i)} - \tilde{U}_{k-1}^{(i)T} \frac{\partial K^{(i)}}{\partial a} \tilde{U}_{k-1}^{(i)} \right|}{\left| \tilde{U}_k^{(i)T} \frac{\partial K^{(i)}}{\partial a} \tilde{U}_k^{(i)} \right|} < \varepsilon^* \tag{39}$$

with

$$\tilde{U}_k^{(i)} = U_{0\ k}^{(i)} + P(dx_k^{(1)} + B_k^{(i)} y_k^{(i)}), \tag{40}$$

where $k$ presents the number of iteration in MG. A suggested value of $\varepsilon^*$ is set to $10^{-2}$, which can be also adjusted according to specific problems [7].

### 3.3.3. Hybrid optimizer

It is important for topology optimization to select a suitable optimization algorithm to obtain reasonable configuration. Sometimes, the non-monotonous behavior may appear during the topology optimization procedure for minimum-compliance problems. This phenomenon results in oscillation especially in the vicinity of the optimal point, which postpones convergence. To handle this problems, monotonous convex approximation MMA and non-monotonous convex approximation GCMMA are integrated as a hybrid optimizer to improve the performance of optimization. More detailed information about these theories can be found in [36, 37]. Both convex approximations based on MMA and GCMMA are presented in Fig. 5.

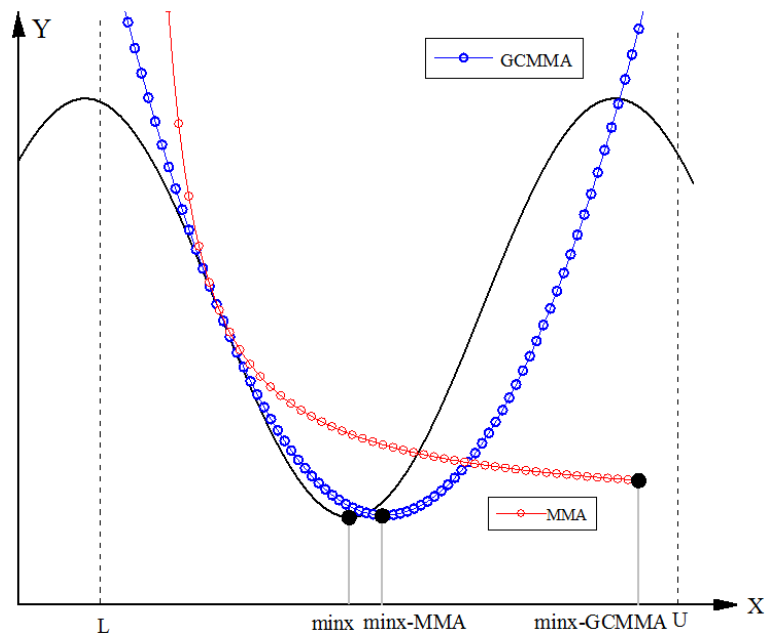

**Fig. 5.** MMA and GCMMA approximate original topology problem monotonously and non-monotonously

Not only the monotonous approximation of MMA accounts for the oscillation near the optimal point, but also the intrinsic character of MMC affects robustness of optimization. Because of the diversity of design variables, the optimization is not easy to converge. Sometimes structural evolutions only generate around a fixed small part of the whole structure recurrently in the neighborhood of optimal point. In addition, the frequent switches between full analysis and reanalysis also intensify the unreliability of optimization. Thus, it is necessary to propose a hybrid algorithm based on MMA and GCMMA [40], to tackle the large-scale and complex topology optimization problems. On the one hand, MMA used in topology optimization can contribute to quick convergence and rapid descent of objective function at the first few iterations. Whereas, it results in numerical oscillation near the optimum which slows down convergence and calculation speed. On the other hand, GCMMA can obtain a stable solution and speed up convergence near the optimal point, but when it is used solely the whole convergence process is very slow. Under such circumstance, the MMA is used in the beginning of topology optimization, primarily when the design point is far from the optimum, and GCMMA plays the main role when the design point is near the optimum. Therefore, it is imperative to state a certain criterion to switch from MMA to GCMMA adaptively during topology optimization. Owing to the focus of avoiding oscillation, the information of it can be used as a switch condition. When the oscillation value of the objective function is within a fixed range, the optimizer will be changed from MMA to GCMMA. The objective values of current step and previous two steps are applied to represent the criterion:

$$-\delta < \frac{f_0(\bm{x}^{(i-2)}) - f_0(\bm{x}^{(i-1)})}{(|f_0(\bm{x}^{(i-2)})| + |f_0(\bm{x}^{(i-1)})|)/2} \times \frac{f_0(\bm{x}^{(i-1)}) - f_0(\bm{x}^{(i)})}{(|f_0(\bm{x}^{(i-1)})| + |f_0(\bm{x}^{(i)})|)/2} < 0, \qquad (41)$$

where $\delta$ is a small positive real number. In this study, it is usually defined as $\delta < 0.002$, which was also obtained by a lot of numerical experiments and could be changed slightly near the defined value.

## 4. Numerical examples

In this section, several benchmark numerical results are presented to demonstrate the effectiveness of the IRA-MMC based topology optimization, with a hybrid optimizer MMA/GCMMA. The discussed problems can be classed into two mainstream categories: minimum compliance problems and compliant mechanism problems. It is shown that accurate results can be achieved efficiently by the proposed approach. All the examples are implemented in MATLAB, and corresponding MATLAB codes are executed in MATLAB R2017a on a desktop computer with Intel Xeon E3-1230 V2 CPU @3.30GHZ and 24GB RAM.

## 4.1. The cantilever example

In this example, the minimum compliance problem for a cantilever beam is addressed by the IRA-MMC based topology optimization versus the original MMC-based topology optimization [13]. As for the geometry information and load constraints, the setup can be seen from Fig. 6, where a concentrated vertical load is enforced at the middle point of the right side, and the left side of the design domain is fixed wholly. The design domain with width $DW = 2$ and height $DH = 1$, is discretized into $80 \times 40$, $100 \times 50$ and $160 \times 80$ FEM meshes in sequence by four-node bi-linear elements. Material properties about Young's modulus and Poisson's ratio are set to $E = 1$ and $v = 0.3$ respectively. The allowed volume fraction is defined as $v = 0.4$. During the implementation, the switch conditions are $\eta = 0.11$, $\varepsilon^* = 10^{-2}$ and $\delta = 0.002$. Both topology optimizations with full analysis and reanalysis adopt the same initial components layout as Fig. 7.

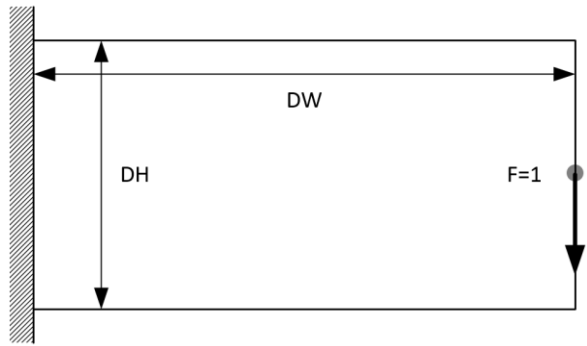

**Fig. 6.** 2D cantilever beam

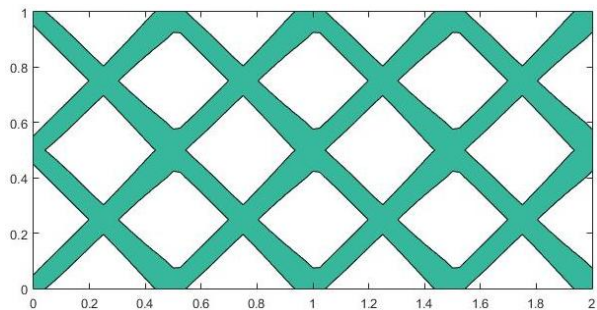

**Fig. 7.** The initial design of the cantilever beam for both full analysis and reanalysis

The final optimization layouts in relation to the proposed approach and the

MMC-based topology optimization, are presented in Fig. 8. It is obvious that structural configurations obtained by the proposed approach are nearly identical to the configurations optimized by the full analysis, and the consistency is not affected by extend scales of grids. Reference to Fig. 9, the merged convergence diagram corresponding to the grid scale of $80\times 40$ can present a clear procedure of optimization. Due to the relatively large oscillation during the first several iterations, the convergence diagram is intercepted to retain the observability of the convergent procedure. The convergence curve in red coincides well with the blue curve, which proves that the proposed method keeps robust and reliable as the original method. Moreover, the convergence curve in red terminates earlier than the blue one, which indicates that the proposed approach can satisfy the stopping criterion more quickly than the MMC-based topology optimization. For demonstrating the efficiency of the proposed approach, the same stopping criterion that requires a maximum change of $10^{-3}$ for all design variables is adopted.

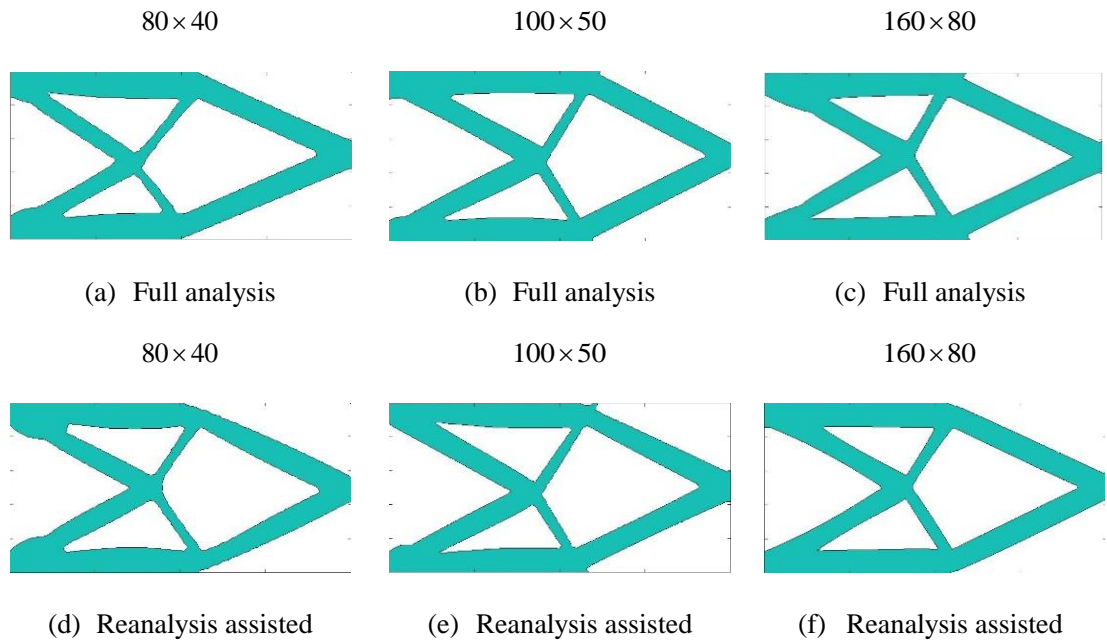

（a） Full analysis     （b） Full analysis     （c） Full analysis

（d） Reanalysis assisted     （e） Reanalysis assisted     （f） Reanalysis assisted

**Fig. 8.** Benchmarked optimized layouts of the cantilever beam between full analysis and IRA-MMC. The scales of FEM elements are refined sequentially from $80\times 40$, to $100\times 50$, and until $160\times 80$

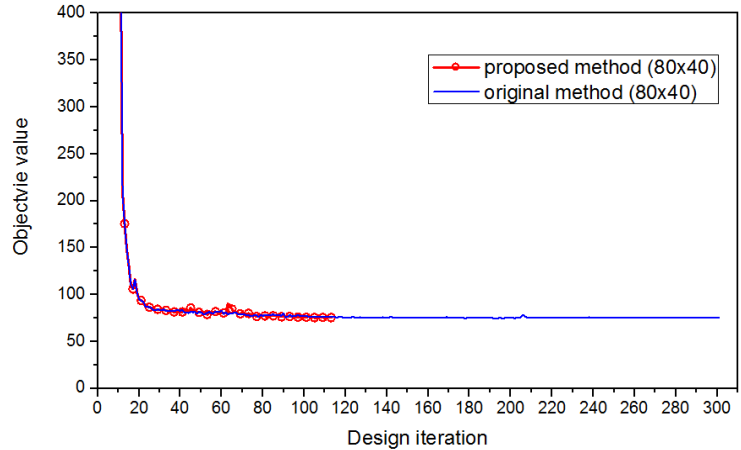

**Fig. 9.** Comprehensive convergence diagrams of objective functions with respect to full analysis (blue) and reanalysis (red) with scale of $80 \times 40$

Table 2 Optimization results of the cantilever beam are shown with various scales of FEM elements. The scale value of each row from top to the bottom corresponds to $80 \times 40$, $100 \times 50$, and $160 \times 80$ respectively

| Approach | Objective | Constraint | Iteration | CPU time (s) |
|---|---|---|---|---|
| IRA-MMC-MMA/GCMMA (three scales) | 75.3749 | $-0.166 \times 10^{-3}$ | 113 | 270.8 |
|  | 74.7304 | $-0.886 \times 10^{-5}$ | 127 | 555.6 |
|  | 75.1425 | $-0.12 \times 10^{-2}$ | 100 | 1008.9 |
| Full analysis (three scales) | 75.3619 | $-0.231 \times 10^{-4}$ | 301 | 777.2 |
|  | 74.5413 | $0.659 \times 10^{-5}$ | 301 | 1028.1 |
|  | 75.1258 | $-0.543 \times 10^{-4}$ | 670 | 4521.5 |
| Percentage difference (three scales) | 0.01725% | - | -62% | -65.16% |
|  | 0.25368% | - | -43.56% | -45.96% |
|  | 0.02223% | - | -85.07% | -77.69% |

Apparently, it can be concluded that consistent optimized structures as full analysis can be obtained by the proposed method with more efficient convergence. Moreover, the sufficient accuracy and improved efficiency of the proposed method can be further validated by specific data, which are presented in Table 2. Remarkable saving in computer time has been obtained with up to 77.69% without reducing accuracy. Although all of the three minimum compliances obtained by the proposed method with gradually refined grids are slightly higher than the full analysis, the biggest percentage difference does not exceed 0.3%. Therefore, the values of objective functions obtained by the proposed approach can be regarded as precise results. Consulting the third column from Table1, the volume constrains in two different methods all achieve the requirement, except for the second one solved by the original method. This can be recognized as a further proof of merits in accuracy and stability about the proposed method.

**4.2. The L-Shape beam example**

Another minimum compliance problem for the L-Shape beam example is also popular in the literature of the topology optimization. The specific geometry of the design domain, and load constraints are depicted in **Fig. 10**. The top of the domain is fixed and a concentrated vertical downward force is enforced at the middle of the right side. The whole domain will be sequentially discretized into 4864, 7600 and 19456 FEM meshes by four-node bi-linear elements. Young's modulus and Poisson's ratio are set to $E=1$ and $v=0.3$ respectively. The allowed volume fraction is defined as $v=0.3$. In this example, the switch conditions are $\eta=0.05$, $\varepsilon^*=10^{-2}$ and $\delta=0.001$. Both topology optimizations with full analysis and reanalysis adopt the same initial components layout as Fig. 11.

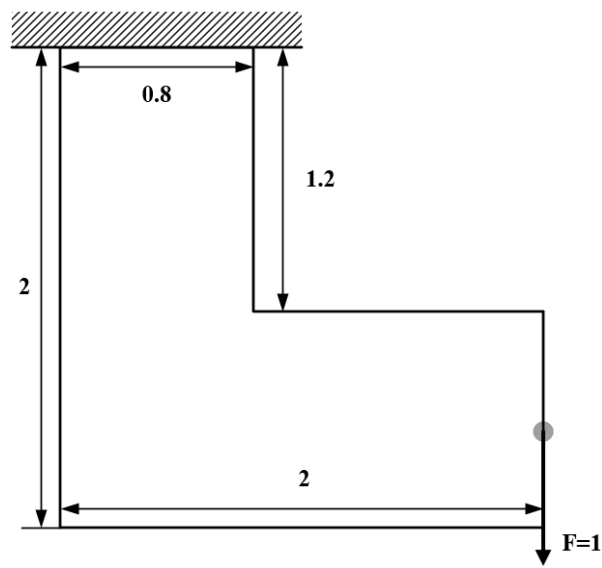

**Fig. 10.** The L-shape beam

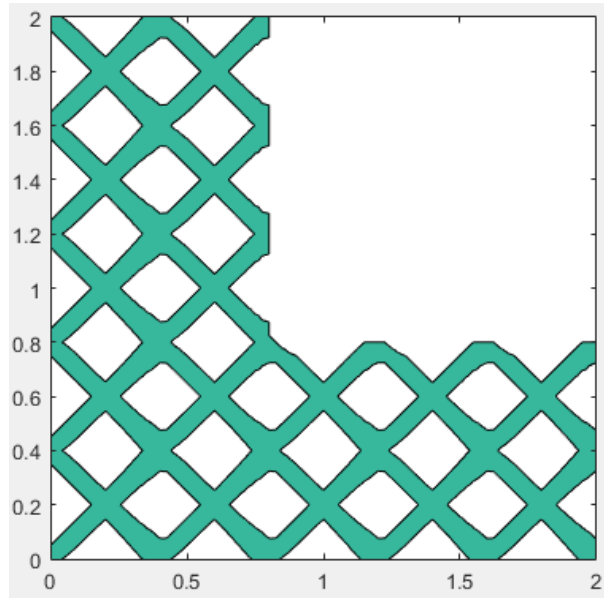

**Fig. 11.** The initial design of the L-shape beam for full analysis and reanalysis

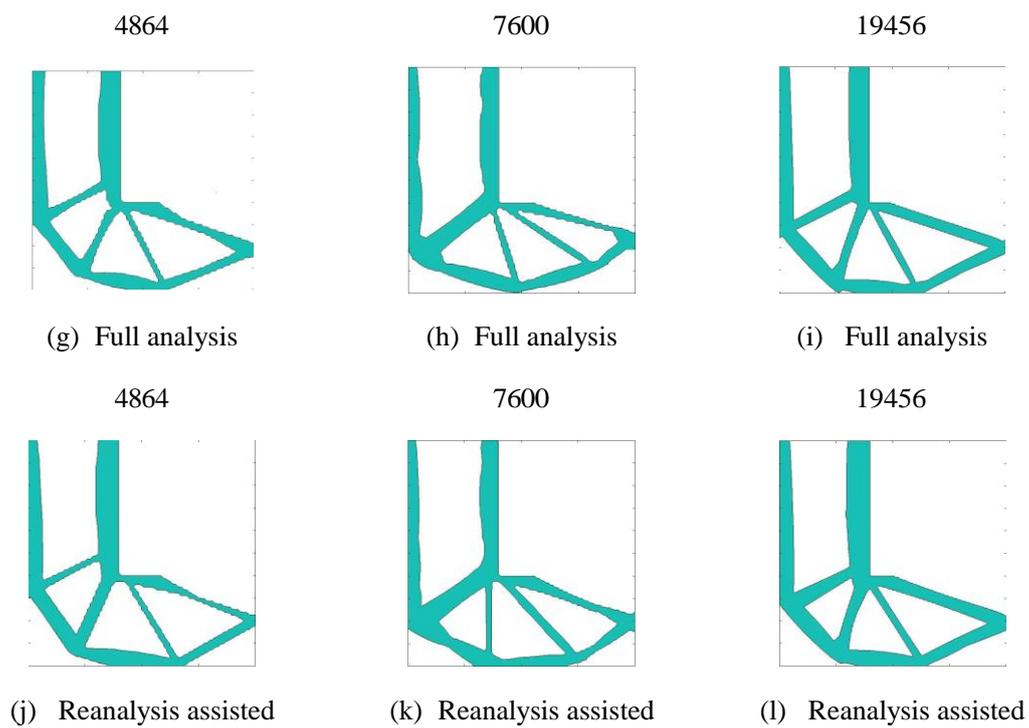

**Fig. 12.** Benchmark optimized layouts of the L-shape beam between full analysis and IRA-MMC. The scales of FEM elements are refined sequentially from 4864, to 7600, and until 19456

According to Fig. 12, all the three optimized structures of different scales in the proposed method are similar to the corresponding structures in the original method. Compared with the original method, optimized structures by the IRA-MMC based topology optimization have more smooth boundaries and more uniform thickness, which are beneficial to manufacturability. Although the convergent procedures in Fig. 13 appear slight oscillations by the proposed method, the accuracy and efficiency are still considerable. From

Table 3, it can be seen that the computational time is saved at least 65.15% among three different scales of models with the maximum objective difference of 5.959%. Despite there are some flaws, it still demonstrates the effectiveness of the IRA-MMC based topology optimization method.

Table 3 Optimization results of the cantilever beam are shown with various scales of FEM elements. The value of each row from top to the bottom corresponds to $80\times80$, $100\times100$, and $160\times160$ respectively

| Approach | Objective | Constraint | Iteration | CPU time (s) |
|---|---|---|---|---|
| IRA-MMC-MMA/GCMMA (three scales) | 239.3844 | 0.0013 | 284 | $2.36\times10^3$ |
| | 234.9183 | $3.57\times10^{-4}$ | 242 | $5.68\times10^3$ |
| | 236.8588 | 0.0024 | 101 | $2.28\times10^3$ |
| Full analysis (three scales) | 238.4589 | $6.53\times10^{-4}$ | 600 | $1.49\times10^4$ |
| | 249.8042 | $1.50\times10^{-4}$ | 600 | $1.63\times10^4$ |
| | 235.3525 | $2.48\times10^{-4}$ | 600 | $1.78\times10^4$ |
| Percentage difference (three scales) | 0.3881% | - | -52.67% | -84.16% |
| | -5.959% | - | -59.67% | -65.15% |
| | 0.6400% | - | -83.17% | -87.19% |

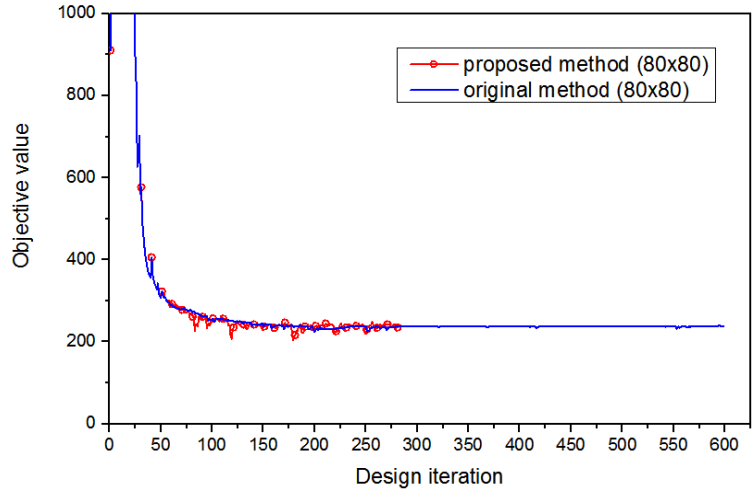

**Fig. 13.** Comprehensive convergence diagrams of objective functions with respect to full analysis (blue) and reanalysis (red) with scale of $80 \times 80$

### 4.3. The compliant mechanism example

The proposed approach is also competent in solving non-self-adjoint problems, just as the compliant mechanism problem. This kind of problems is more complex and sensitive to parameters setting of MMA and switch criterions. The problem setup is shown in Fig. 14. The compliant mechanism example Constants of springs at the input/output points are $k_{in} = k_{out} = 0.1$ respectively. For a qualitative examination, the whole design domain is also discretized as the first examples with $80 \times 40$, $100 \times 50$ and $160 \times 80$ FEM meshes by four-node bi-linear elements. The Young's modulus $E$ and Poisson's ratio $v$ are set to 1 and 0.3 respectively. The volume fraction is set to $v = 0.3$. Parameters of switch criterions can be set as $\eta = 0.13$, $\varepsilon^* = 10^{-2}$ and $\delta = 0.002$. Initial layout for optimization can be consulted from Fig. 15.

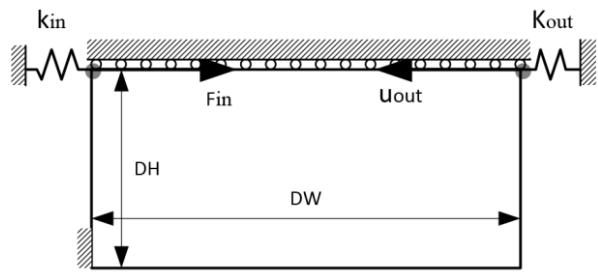

**Fig. 14.** The compliant mechanism example

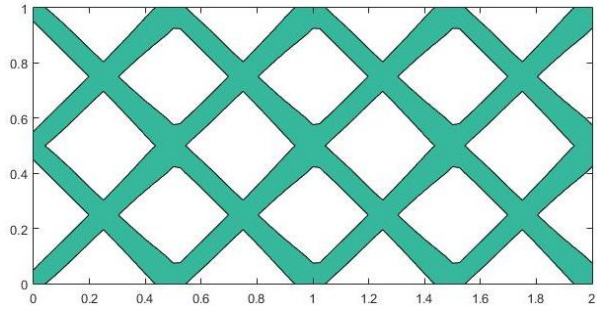

**Fig. 15.** The initial design for the compliant mechanism example

The IRA-MMC based topology optimization successfully circumvents a well-known strong local minimum at the objective value of zero about compliant mechanism problems. All the three scales of compliant mechanism examples achieve convergence at the end, and acquire reasonable layouts. Reference to Fig. 16, along with mesh refinement the proposed method obtains nearly consistent layouts after optimization. This phenomenon is as same as the final layouts obtained by the MMC-based topology optimization. This demonstrates the proposed method inherits stability of the original method. Although the layouts optimized by the proposed method are slightly different from the corresponding layouts by the original approach, they share almost the identical objective values with the maximal bias of 1.22% and the minimal bias of 0.14%. It can be indicated from Fig. 17 that the red curve coincides with the blue one well and satisfies convergence criteria earlier.

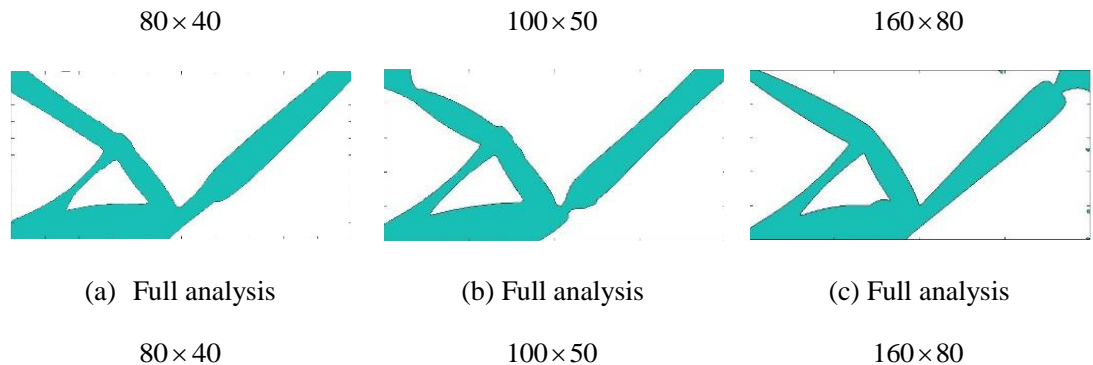

$80 \times 40$        $100 \times 50$        $160 \times 80$

(a) Full analysis     (b) Full analysis     (c) Full analysis

$80 \times 40$        $100 \times 50$        $160 \times 80$

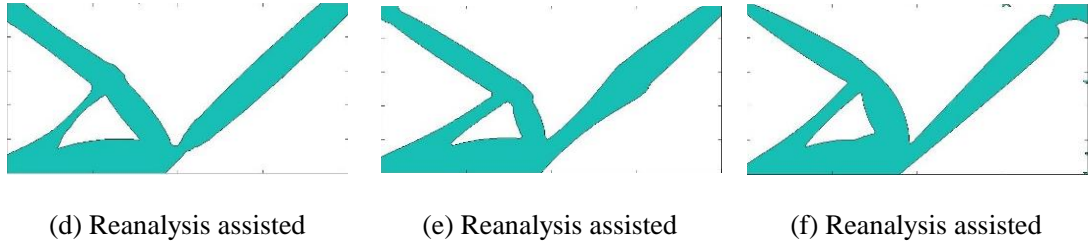

(d) Reanalysis assisted      (e) Reanalysis assisted      (f) Reanalysis assisted

**Fig. 16.** Benchmark optimized layouts of the compliant mechanism example between full analysis and IRA-MMC. The scales of FEM elements are refined sequentially from $80 \times 40$, to $100 \times 50$, and until $160 \times 80$

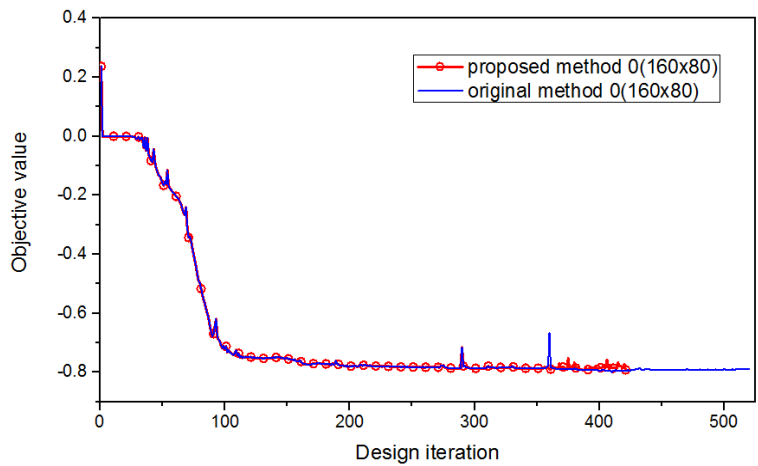

**Fig. 17.** Comprehensive convergence diagrams of objective functions with respect to full analysis (blue) and reanalysis (red) with scales of $160 \times 80$

Besides the appearances of objects, specific data can validate the efficiency and accuracy of the proposed method more explicitly. Under the condition of satisfying volume constraint, the proposed approach achieves nearly exact objective values with only trivial bias, which can be recognized through the second column of Table 4. However, the constraint values solved by original method are greater than zero when it comes to the scale of $100 \times 50$, and $160 \times 80$. This further demonstrates that the proposed method can obtain more reliable results. Moreover, all the three examples in different scales successfully save computational cost, which acquires at least 50% reduction. The effective performances of the proposed approach prove that it is also appropriate for non-self-adjoint problems.

Table 4 Optimization results of the compliant mechanism are shown with various scales of FEM elements. The values of each row from top to the bottom correspond to $80\times 40$, $100\times 50$, and $160\times 80$ respectively

| Approach | Objective | Constraint | Iteration | CPU time (s) |
|---|---|---|---|---|
| IRA-MMC-MMA/GCMMA (three scales) | -0.7873 | $-7.1359\times 10^{-5}$ | 350 | $1.2470\times 10^{3}$ |
|  | -0.7932 | $-0.0015$ | 305 | $1.4862\times 10^{3}$ |
|  | -0.7919 | $-0.0012$ | 421 | $4.1917\times 10^{3}$ |
| Full analysis (three scales) | -0.7778 | $-3.1668\times 10^{-4}$ | 612 | $4.1598\times 10^{3}$ |
|  | -0.7991 | $3.2966\times 10^{-5}$ | 717 | $5.7725\times 10^{3}$ |
|  | -0.7908 | $4.9925\times 10^{-6}$ | 520 | $8.2523\times 10^{3}$ |
| Percentage difference (three scales) | 1.22% | - | -42.81% | -70.02% |
|  | -0.74% | - | -57.46% | -74.25% |
|  | 0.14% | - | -19.04% | -49.21% |

## 5. Conclusions

In this study, an iterative reanalysis approximation method, IRA, is integrated into structural topology optimization based on the MMC for the static analysis. A hybrid optimizer constructed by MMA and GCMMA substitutes for the single MMA. In comparison with the original method, the IRA-MMC based topology optimization reduces the number of matrix factorizations, decreases dimensionalities of equilibrium equations directly solved and accelerates convergence without compromising accuracy. Suitable switching criteria for analysis solvers and optimizers play a critical part in improving efficiency and inheriting robustness.

The proposed approach is in contrary to the original one by three benchmark problems on a sequence of increasingly refined meshes. The largest scale of problem has more than 26 thousand DOFs. In spite of the high efficiency provided by direct sparse solvers in 2-D problems, the proposed approach still runs more than twice faster compared with the original method, especially in the first and the last example.

In addition, all the three examples obtain stable and consistent design structures as the mesh refinement. The remarkable improvement on computational cost and explicit design structures demonstrate the proposed approach accommodates the minimum compliance problems and compliant mechanism problems. Furthermore, it successfully improves efficiency with reliable procedure and inherits the explicitness of the MMC-based topology optimization simultaneously.

# Acknowledgments

This work has been supported by Project of the Key Program of National Natural Science Foundation of China under the Grant Numbers 11572120.